# Arrow graphs and their focal curves

Rogier Bos & Filip Cools

## Abstract

*Recent research in mathematics education has renewed interest in using arrow graphs to represent functions in secondary education. An arrow graph depicts pairs $(x, f(x))$ as arrows connecting two parallel axes. When examining such graphs, a striking curve appears, formed by the envelope of these arrows. This curve, which we named the focal curve, is the central object of study in this paper. We investigate focal curves in depth by deriving parametrizations and equations for various common functions and exploring their geometric properties. Many of these properties can be explained by viewing the focal curve as the dual of the function's graph when interpreted in the projective plane. This leads us to examine focal curves from the perspective of projective geometry. Finally, we analyze how focal curves change under transformations and compositions of the original functions.*

## 1. Introduction

The most common visualization of a function is undoubtedly its graph. However, more obscure options are available, such as the arrow graph (Richmond, 1963; Brieske, 1978), also known as the arrow diagram, mapping diagram, or nomogram. The arrow graph of a function $f: \mathbb{R} \to \mathbb{R}$ is drawn in a parallel axes system, consisting of two parallel axes—the input and the output axes. It consists of a set of arrows $x \to f(x)$ from the values $x$ on the input axis to the corresponding value $f(x)$ on the output axis, where the input values are usually chosen at regular intervals. See some examples of arrow graphs in Figure 1. Recent studies explored the application of arrow graphs in secondary education (Bos & Brinks, 2024; Wei et al., 2024, 2025). The vectors in the arrow graphs in Figure 1 envelop a mysterious curve—a hyperbola on the left and a circle on the right—which plays the central role in this article. There is a beautiful relationship between this curve, which we baptized the *focal curve*, and the derivative of the function (Bridger, 1996). We explore focal curves in detail in this paper. Indeed, we explain how it is parametrized; compute various examples; show how it relates to the dual curve of graph; and, finally, discuss how it transforms as the function transforms.

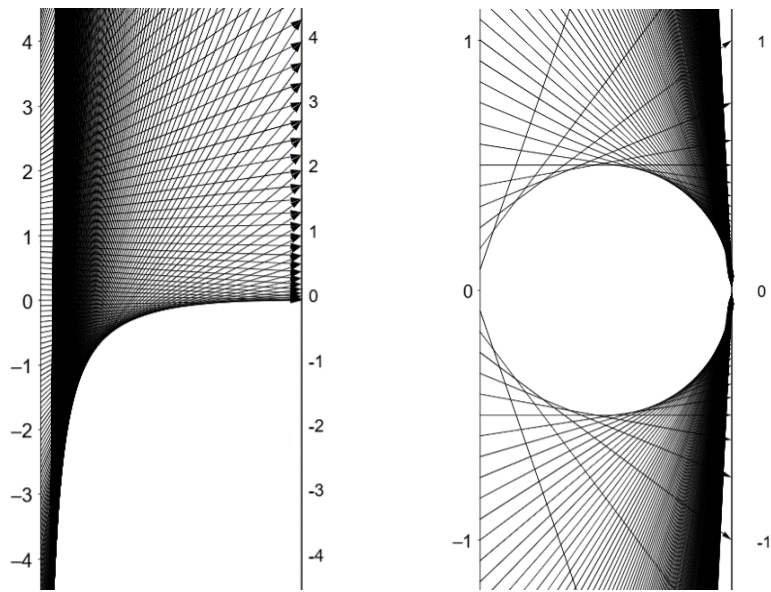

Figure 1. The arrow graph of the function given by $f(x) = x^2$ and $g(x) = \frac{1}{4x}$.

Exploring the arrow graphs of other functions in this GeoGebra app is interesting and fun: https://www.geogebra.org/m/gh7zcw93. One observation you cannot miss is that the grey lines in the arrow graph of a *linear* function intersect at one point (unless the gradient equals 1, in which case the vectors are parallel). Indeed, for a linear function $f(x) = a x + b$, any interval of width $\Delta x$ on the input axis is enlarged to an interval on the output axis of width $\Delta y = a \Delta x$ with respect to the intersection point $F$ (see Figure 2). This point is called the focus of the function, and its position uniquely determines the linear function. As mentioned, the parameter $a$ is the enlargement factor and $b$ is the value where the line through the focus and the zero in the input axis intersects the output axis.

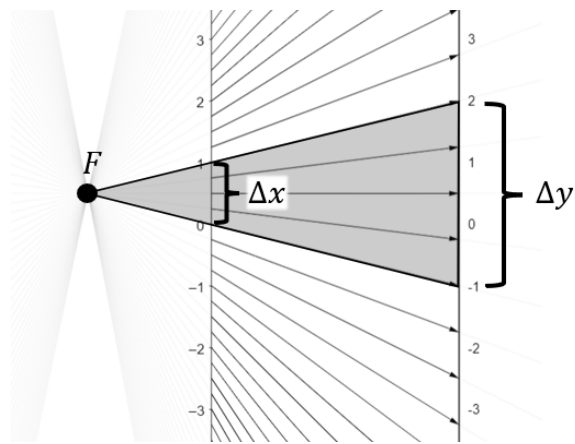

Figure 2. Arrow graph for the linear function $f(x) = 3x - 1$. The gradient 3 is the enlargement factor of the corresponding intervals with respect to the focal point $F$.

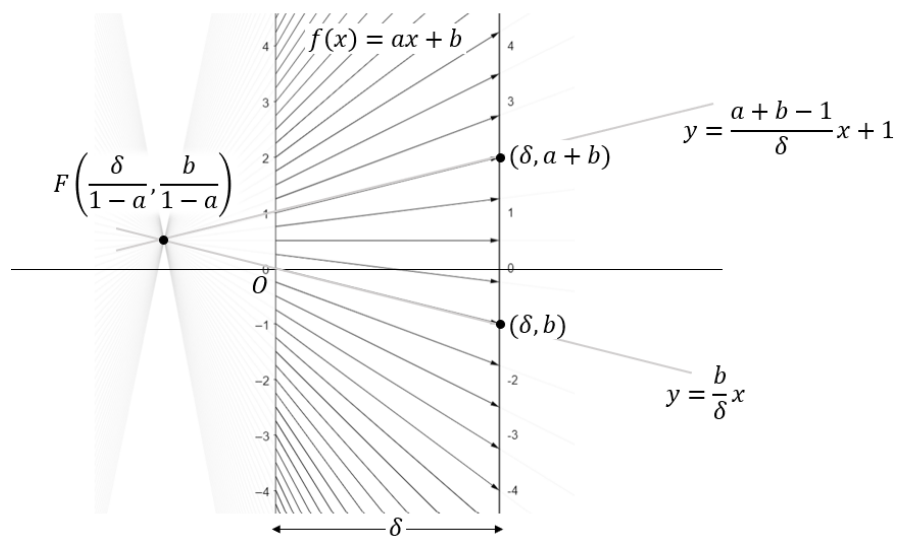

Figure 3. The arrow graph placed in an Cartesian axes system. Two well-chosen lines help to compute the Cartesian coordinates of the focus in terms of the parameters of the linear function $f$.

One can place the parallel axes system in an orthogonal axes system such that the input axis coincides with the vertical axis, and the output axis is placed on the vertical line associated with $x = \delta$ for a positive real number $\delta$. We set $\delta = 1$, unless stated otherwise. In this GeoGebra app, you can explore the relation between the focus coordinates and the linear function parameters: https://www.geogebra.org/m/gq5m8dsd. By intersecting the two well-chosen lines a straightforward computation gives the coordinates of the focus in terms of the parameters $a, b$ and $\delta$: $F\left(\frac{\delta}{1-a}, \frac{b}{1-a}\right)$, see Figure 3 (Bridger,1996). The Cartesian plane $\mathbb{R}^2$ can be embedded in the projective plane $\mathbb{P}^2$ by mapping $(x, y)$ to $(x : y : 1)$. Under this embedding the focal point corresponds to the point with projective coordinates $(1 : b : 1 - a)$, choosing $\delta = 1$. The functions with gradient $a = 1$ with parallel vectors have a focal point at infinity, corresponding to points $(1 : b : 0)$ in the projective plane.

## 2. the focal curve

In this section, we discuss the relation between the enveloping curve emerging from the arrow graph of a function $f$ and its derivative $f'$. A function is differentiable in $x_0$, whenever the graph is locally linear near the point $(x_0, f(x_0))$. In the arrow graph, this means that on an interval around $x_0$ on the input axis, the arrow graph looks like the arrow graph of a linear function, that is, the lines through the vectors on that interval go approximately through one point, see Figure 4. This point we call the local focus for $x_0$.

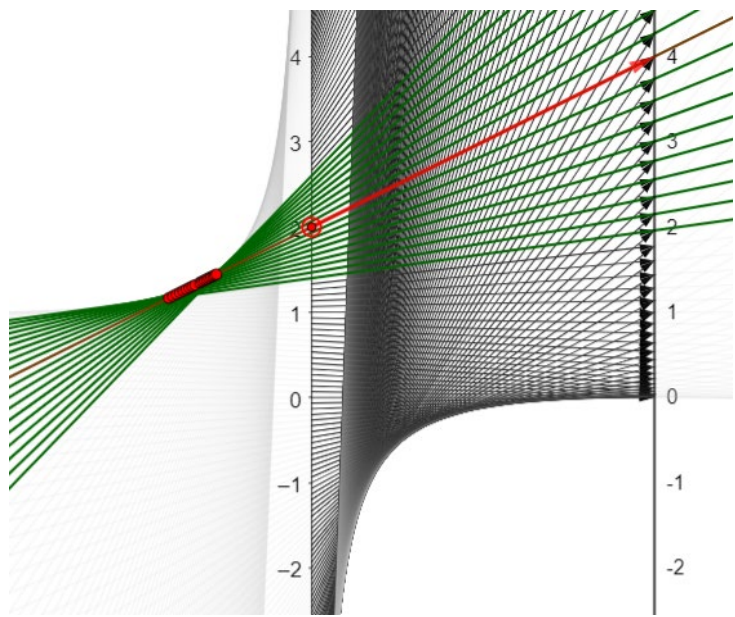

Figure 4. The lines through the arrows on an interval near $x_0 = 2$ go approximately through one point: the local focus.

The best linear approximation of the function $f$ near $x_0$ is given by $y = f(x_0) + f'(x_0)(x - x_0)$. The focus for this line is the local focus $F_{x_0}$ for the curve at $x_0$ and, using the coordinates found in Figure 3, we find the coordinates $\left(\frac{1}{1-f'(x_0)}, \frac{f(x_0) - x_0 f'(x_0)}{1-f'(x_0)}\right)$ for $F_{x_0}$ (Bridger, 1996). Together, the local foci constitute a curve called the *focal curve*.

**Proposition:** The enveloping curve of the lines through the vectors of an arrow graph equals the focal curve.

**Proof:** The only thing left to prove is that the vectors from $(0, x)$ to $(1, f(x))$ are tangent to the focal curve parametrized by $(x(t), y(t)) = \left(\frac{1}{1-f'(t)}, \frac{f(t) - t f'(t)}{1-f'(t)}\right)$ for all $x$. Indeed, a small computation shows that

$$(x'(t), y'(t)) = \frac{f''(t)}{\left(1 - f'(t)\right)^2} (1, f(t) - t),$$

which is a scalar multiple of the vector from $(0, t)$ to $(1, f(t))$, which equals $(1, f(t) - t)$.

**Remark 1:** The concept of local focus allows us to interpret $f'(x_0)$ as the enlargement factor for enlargements of intervals on the input axis to the output axis with respect $F_{x_0}$. The notion of enlargement is close to the interpretation of the derivative as the "sensitivity" of the output for variation in the input. The value of $f'(x_0)$ can be determined from the "horizontal" position of the local focus with respect to the vertical axes. Indeed, $f'(x_0) = \Delta o / \Delta i$, see Figure 5 for an example.

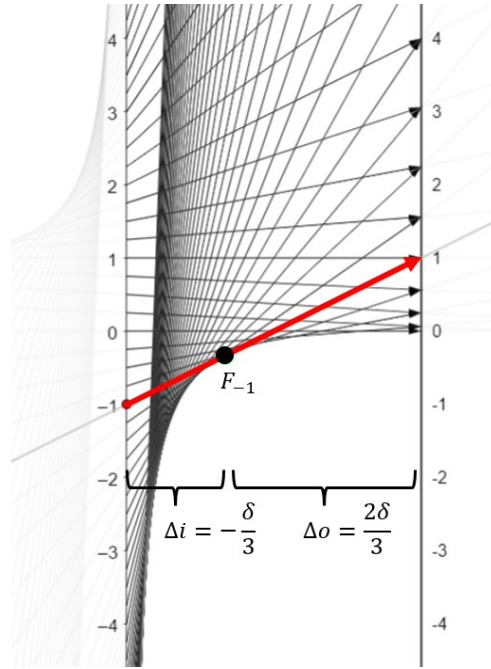

Figure 5. In the arrow graph of $f(x) = x^2$, the local focus for $x = -1$ is found by intersecting the vector $-1 \to 1$ with the focal curve. Then, the derivative $f'(-1)$ is determined as $\frac{\Delta o}{\Delta i} = -2$.

## 3. Some examples of focal curves

The parametrization $(x(t), y(t)) = \left(\frac{1}{1-f'(t)}, \frac{f(t)-tf'(t)}{1-f'(t)}\right)$ of the focal curve of a function $f$ can be used to compute the focal curves of specific functions. For basic functions, this computation can be done by hand. For more complicated functions, we use the computer algebra system SageMath (The Sage Developers, 2025).

*Example 1 (inverse function, see Figure 1).* Using $f'(t) = -\frac{1}{4t^2}$, we find that $x(t) = \frac{1}{1+\frac{1}{4t^2}}$ and $y = \frac{\frac{1}{2t}}{1+\frac{1}{4t^2}}$. From the former, we derive that $\frac{1}{4t^2} = \frac{1}{x} - 1$. Hence $y^2 = x^2 \left(\frac{1}{x} - 1\right)$, which can be rewritten as $\left(x - \frac{1}{2}\right)^2 + y^2 = \frac{1}{4}$, the equation for a circle of radius 1/2 with center (0,1/2)

*Example 2 (the exponential function).* The focal curve of the function $f(x) = e^x$ is given by the equation $y = (x - 1)\left(1 - \ln\left(\frac{x-1}{x}\right)\right)$. Indeed, using the parametrization, we have that $x(t) = \frac{1}{1-e^t}$, hence $t = \ln\left(\frac{x-1}{x}\right)$. The equation now follows from $y(t) = \frac{(1-t)e^t}{1-e^t} = (1-t)e^t x(t)$.

*Example 3 (the square root function).* The focal curve of the function $f(x) = \sqrt{x}$ has equation $y = \frac{x^2}{4(x-1)}$, with $x < 0$ or $x > 1$. Note that the focal curve is a part of the plane conic given by the equation $4(x-1)y = x^2$, see Figure 6 (right), where however the part for $x < 0$ is not visible.

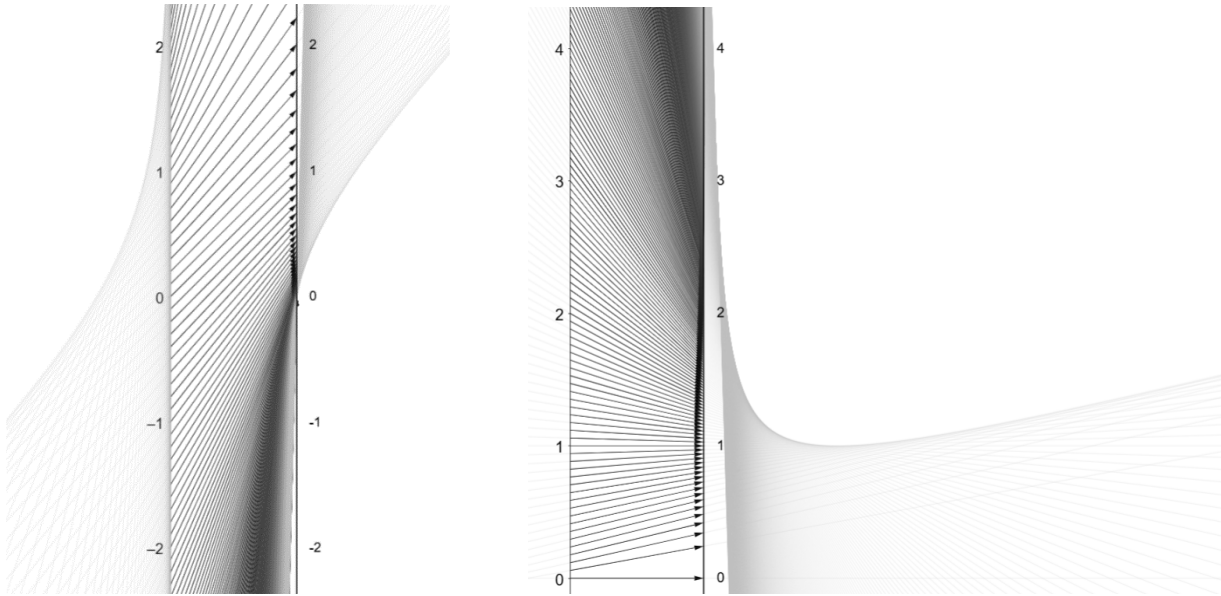

Figure 6. Arrow graphs for the exponential function (left) and the square root function (right), with $\delta = 1$.

*Example 4 (a family of rational functions).* Consider a rational function of the form $f(x) = \frac{ax^2+bx+c}{dx+e}$. Using SageMath, we can eliminate the parameter $t$ from the parameterization of the focal curve. One obtains the equation $(b^2 - 4ac + 4cd - 2be + e^2)x^2 - 2(bd - 2ae + de)xy + d^2y^2 - 2(2cd - be + e^2)x + 2dey + e^2 = 0$, so the focal curve is a plane conic. Moreover, the focal curve intersects the vertical axis $x = 0$ only at the point $\left(0, -\frac{e}{d}\right)$. Hence, the point on the vertical axis of which the $y$-coordinate corresponds to the pole of the original function $f$ (the point at infinity if $d = 0$) belongs to the focal curve and at this point the tangent line is the vertical axis. By looking at some specific cases, we see that different types of conics occur.

- For the function $f(x) = x^2$, the focal curve is the hyperbola given by $(x-1)^2 + 4xy = 0$.
- For the function $f(x) = \frac{1}{4x}$, the focal curve is the circle with center $(1/2, 0)$ and radius $1/2$ (as we have seen in Example 1).
- For the function $f(x) = x + \frac{1}{x}$, the focal curve is the parabola given by $y^2 - 4x = 0$.

Using SageMath, one can also check that each conic which is tangent at the vertical axis occurs as the focal curve of a rational function of the studied form.

## 4. A projective viewpoint: the relation with dual curves

We have seen that the focal curve of a linear function is a point and that the focal curve of the functions in Example 4 (whose graphs are conics) are conics. This suggests that there might be a relation to the theory of dual curves in projective geometry. Recall that the dual curve $C^*$ of a projective plane curve $C$ consists of the projective points $(A : B : C) \in \mathbb{P}^2$ for which $AX + BY + CZ = 0$ is a tangent line to $C$. If we view the graph of the function $f$ as a projective curve $C$, we can consider both its focal curve $F$ and its dual curve $C^*$.

Recall that a line in $\mathbb{R}^2$, corresponding to an equation $y = a\,x + b$, embeds in $\mathbb{P}^2$, through $(x, y) \mapsto (x : y : 1)$, as the projective line defined by $a\,X - Y + b\,Z = 0$. Hence, the dual of such a line corresponds to the point $(a : -1 : b)$. This is not exactly the focus, which embeds in $\mathbb{P}^2$ as $(1 : b : 1 - a)$. However, this leads to the more general insight:

**Proposition:** For any function $f$, the focal curve $F$ and the dual curve $C^*$ are equal up to a projective coordinate transformation of $\mathbb{P}^2$.

**Proof:** The tangent line to the graph of $f$ at a point $(t, f(t))$ is given by the equation $y - f(t) = f'(t)(x - t)$, or in projective coordinates by $f'(t)X - Y + (f(t) - tf'(t))Z = 0$. Hence the dual curve $C^*$ is parametrized by $(f'(t) : -1 : f(t) - tf'(t)) \in C^*$. The focal curve $F$ is parametrized by

$$\left(\frac{1}{1-f'(t)} : \frac{f(t)-tf'(t)}{1-f'(t)} : 1\right) = (1 : f(t) - tf'(t) : 1 - f'(t))$$

Hence, the coordinate transformation $(X : Y : Z) \mapsto (-Y : Z : -X - Y) = (Y : -Z : X + Y)$ maps the dual curve $C^*$ to the focal curve $F$.

**Corollary:** If the function $f$ is rational, then its focal curve is an algebraic curve of degree equal to the class of $C$ (i.e. the number of lines through any fixed point that are tangent to $C$). Indeed, the curve $C$ is algebraic and the dual curve $C^*$ of an algebraic curve is also algebraic of degree equal to the class of $C$.

**Remark 2:** Applied to the rational functions from Example 3, the proposition gives another and more geometric explanation for the fact that the focal curve is a conic: the graph of the function is a conic (defined by the equation $y(d\,x + e) = a\,x^2 + b\,x + c$) and it is known that the dual curve of a conic is a conic, hence also the focal curve is a conic.

*Example 5 (polynomial functions).* If $f$ is a polynomial function of degree $d > 2$, then the projective curve $C$ is singular, and the only singular point $P = (0 : 1 : 0)$ has order $d - 1$. By looking at the degrees of the coordinate functions of the parameterizations in the proof of the Proposition, we can conclude that the dual curve $C^*$ and the focal curve $F$ are plane curves of the same degree $d$. Another way to find the degree is by computing the class of the curve $C$. Therefore, let's fix a point $Q = (a : -1 : 0) \in \mathbb{P}^2$ at infinity. The line at infinity goes through the point $Q$ and is tangent to the curve, since the singular point $P$ also belongs to this line. All other such lines are tangent to the graph of the function $f$ and have a slope equal to $a$, so they correspond to points $(t, f(t))$ with $f'(t) = a$. The latter equation has $d - 1$ solutions, so in total there are $d$ lines that satisfy the condition.

*Example 6 (rational functions).* If $f$ is a rational function with $f(x) = \frac{P(x)}{Q(x)}$ for polynomials $P, Q,$ then the focal curve $F$ is projectively parametrized by

$$(Q(t)^2 : P(t)Q(t) - t(P'(t)Q(t) - P(t)Q'(t)) : P(t)Q'(t) - P'(t)Q(t) + Q(t)^2).$$

The highest degree $d$ of the polynomial coordinate functions in the parametrization depends on the degrees of the polynomials $P$ and $Q$: if $\deg(P) = p$ and $\deg(Q) = q$, then

$$d = \begin{cases} 2q & \text{if } p \leq q + 1 \\ p + q & \text{if } p > q + 1 \end{cases}$$

Hence, $\deg(F) \leq d$ and generically the equality holds (e.g. the polynomials $P$ and $Q$ are at least coprime).

## 5. Focal curves of transformed functions

Previously, through computation, we established the focal curves of some elementary functions. Below, we show how to find the shape of more intricate functions through transformations.

**Proposition**: Suppose the focal curve of a function $g$ is given by the equation $G(x, y) = 0$.

- If $f(x) = g(x) + c$, then its focal curve is given by $G(x, y - cx) = 0$.
- If $f(x) = c\, g(x)$, then its focal curve is given by $G\left(\frac{cx}{1+(c-1)x}, \frac{y}{1+(c-1)x}\right) = 0$.
- If $f(x) = g(x - c)$, then its focal curve is given by $G(x, y + c(x - 1)) = 0$.
- If $f(x) = g(c\, x)$, then its focal curve is given by $G\left(\frac{cx}{1+(c-1)x}, \frac{cy}{1+(c-1)x}\right) = 0$

**Proof:** Let us prove the last statement. A parametrization of the focal curve of $f(x) = g(c\, x)$ is given by

$$(x(t), y(t)) = \left(\frac{1}{1 - g(c\, t)\, c}, \frac{g(c\, t) - t\, g'(c\, t)\, c}{1 - g(c\, t)\, c}\right)$$

From which follows

$$G\left(\frac{c\, x}{1 - (1 - c)x}, \frac{c\, y}{1 - (1 - c)x}\right) = G\left(\frac{c}{1/x - (1 - c)}, \frac{c\, y/x}{1/x - (1 - c)}\right)$$
$$= G\left(\frac{c}{c - g(c\, t)\, c}, \frac{c\, (g(c\, t) - t\, g'(c\, t)\, c)}{c - g(c\, t)\, c}\right)$$
$$= G\left(\frac{1}{1 - g(c\, t)}, \frac{g(c\, t) - c\, t\, g'(c\, t)}{1 - g(c\, t)}\right) = 0.$$

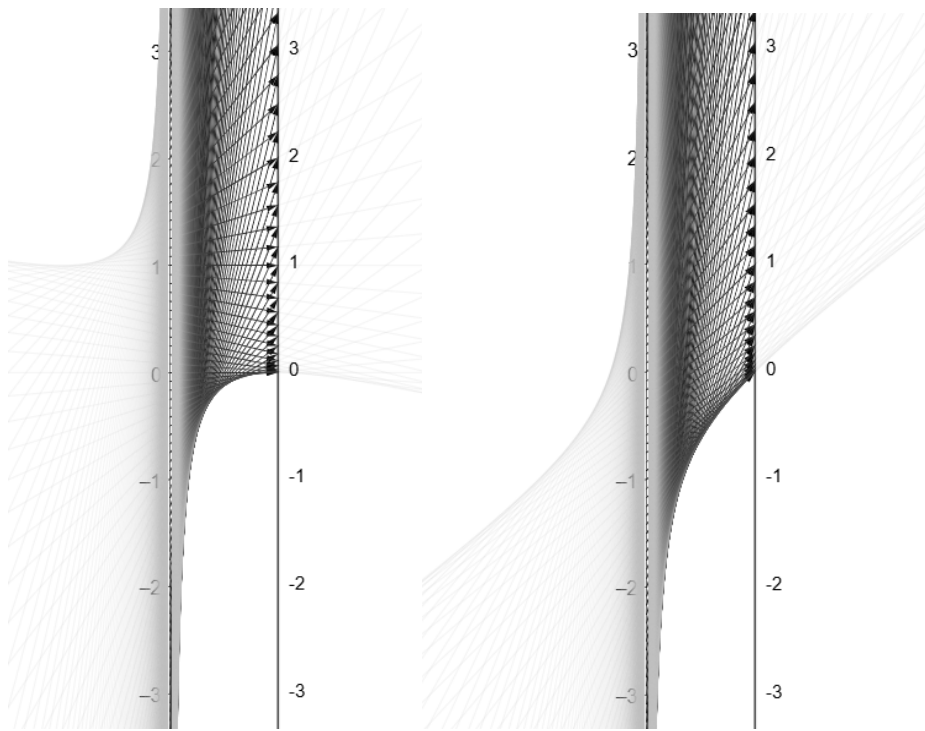

Figure 7. Arrow graphs for $x \mapsto x^2$ and for $x \mapsto (x + 1)^2$. Note how the hyperbolic focal curve transforms according to a shear mapping.

**Remark 3:** If $f$ is a periodic function with period $c$, then its focal curve is invariant under the shear transformation $(x, y) \mapsto (x, y + c(x - 1))$. This transformation fixes the points on the vertical line $x = 1$ and moves all other points vertically by a distance that only depends on the signed distance to the line $x = 1$.

*Example 8 (the sine function).* We consider the periodic function $f(x) = \sin x$ with period $c = 2\pi$ and take a look at its arrow graph and focal curve, see Figure 8. The part of the focal curve that corresponds to the restriction of the function to the interval $[0, 2\pi]$ is colored in blue. It consists of two branches that start off at the point $(0.5, \pi/2)$. In general, each cusp of the focal curve corresponds to an inflection point of the function (this follows from the relation with the dual curve, but can also be checked using the parametrization, and its derivative). In this case, the point $(0.5, \pi/2)$ of the focal curve corresponds to the inflection point $(\pi, 0)$ of the function $f$.

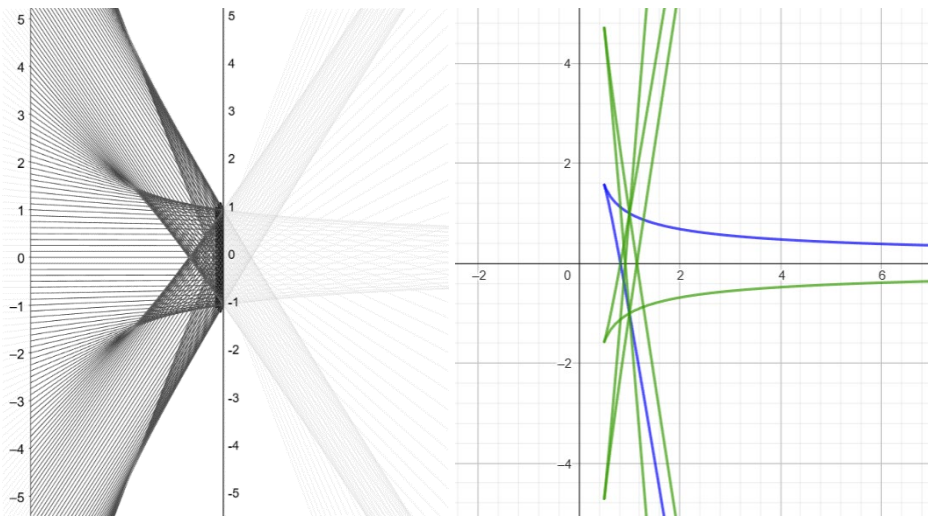

Figure 8. Close-up of arrow graph for $x \mapsto \sin(x)$ (left) and a part of its focal curve (right). In blue is the part that corresponds to the interval $[0, 2\pi]$ as a subset of the domain of the function; other parts are colored in green.

## 6. compositions of functions

A nice feature of arrow graphs is the visualization of the composition of functions by juxtaposing two parallel axes systems, see Figure 9.

**Proposition.** Suppose $f$ and $g$ are linear functions. Then the focal points $F_f$, $F_g$ and $F_{g \circ f}$ of respectively $f$, $g$ and the composition $g \circ f$ are colinear. Moreover, if $f(x) = a x + b$ and $g(x) = c x + d$, then $(1 - t) F_f + t F_g = F_{g \circ f}$, with $t = \frac{1-c}{1-ac}$.

**Proof:** Consider the line through the foci $F_f$ and $F_g$. There are vectors exactly on this line for both the arrow graph of $f$ and $g$. The addition of these vectors gives a vector on the arrow graph of $g \circ f$. The line through this vector goes through $F_{g \circ f}$ equals line $F_f F_g$, which proves the first point. The second point is proved by direct computation, using $F_f = \left(\frac{1}{1-a}, \frac{b}{1-a}\right)$, $F_g = \left(\frac{1}{1-c} + 1, \frac{d}{1-c}\right)$ and $F_{g \circ f} = \left(\frac{2}{1-ac}, \frac{bc+d}{1-ac}\right)$. Note that the parallel axes system of $g$ is translated by 1 to the right, hence the "+1" in the first coordinate of $F_g$, and the parallel axes system for $g \circ f$ has width $\delta = 2$, hence the 2 in the first coordinate of $F_{g \circ f}$.

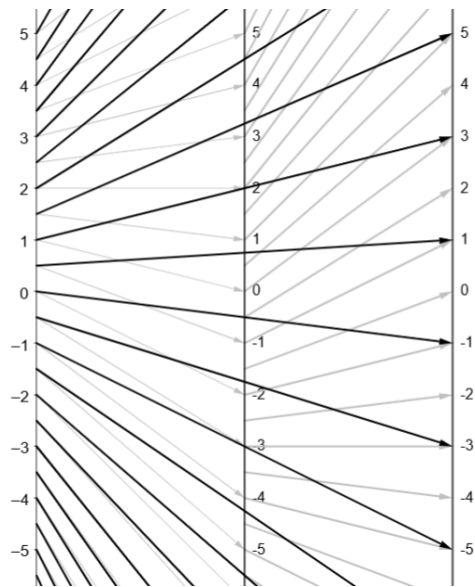

Figure 9. The composition $g \circ f$ of the functions $f(x) = 2x - 2$ and $g(x) = 2x + 3$

In summary, we have seen how a function can be presented graphically by its arrow graph. In the arrow graph, the local foci reveal the derivative function as an enlargement factor. Together, the local foci form the focal curve, a curve related to the dual curve by a coordinate transformation. Focal curves can be computed explicitly and show properties of the function and its derivative in an unusual and interesting way. We believe this allows an enriched understanding of functions, which deserves further exploration, particularly in educational settings.

Rogier Bos
Freudenthal Institute
Utrecht University
Princetonplein 5, Utrecht, the Netherlands
r.d.bos@uu.nl

Filip Cools
Katholiek Onderwijs Vlaanderen
Guimardstraat 1, Brussels, Belgium
filip.cools@katholiekonderwijs.vlaanderen